\documentclass[twoside,10pt]{article}
\usepackage[]{amsmath,amsthm,amssymb}
\usepackage[latin1]{inputenc}

\begin{document}

\begin{center}{{\Large\bf  A new Kontorovich-Lebedev-like transformation} }\end{center}
\vspace{0,5cm}
\begin{center}{Semyon  YAKUBOVICH}\end{center}

\markboth{\rm \centerline{ Semyon  YAKUBOVICH}}{}
\markright{\rm \centerline{Kontorovich-Lebedev-Like Transformation}}

\begin{abstract} {\noindent A different application of the familiar integral representation
for the modified Bessel function drives to a new
Kontorovich-Lebedev-like integral transformation of a general
complex index. Mapping and operational  properties, a convolution
operator  and inversion formula are established. Solvability
conditions and explicit solutions of the corresponding class of
convolution integral equations are exhibited.}

\end{abstract}
\vspace{4mm}
{\bf Keywords}: {\it    Kontorovich-Lebedev transform, modified
Bessel functions, Mellin transform, Laplace transform, convolution,
integral equations of the convolution type}

\vspace{2mm}

 {\bf AMS subject classification}:  44A15, 33C05, 33C10, 33C15

\vspace{4mm}

\section {Introduction}

As it is known [2], Vol. II, the modified Bessel function
$K_z(2\sqrt x)$ can be represented by the following integral
$$K_z(2\sqrt x)= {x^{-z/2}\over 2}\int_0^\infty e^{-t - {x\over
t}}t^{z-1} dt, \ x > 0,\eqno(1.1)$$
where $z=\nu+i\tau$ is a complex number. As it is easily seen,
integral (1.1) converges absolutely for any $x \in \mathbb{R}_+, z
\in \mathbb{C}$ and represents an entire function by $z$. Formula
(1.1) can be written with the use of the Parceval relation for the
Mellin transform [6], which leads to the integral representation
$$ 2 x^{z/2}K_z(2\sqrt x)= {1\over 2\pi i}
\int_{\gamma-i\infty}^{\gamma+i\infty} \Gamma(s+z)\Gamma(s) x^{-s}
ds,\ x > 0, \eqno(1.2)$$
where $\Gamma(w)$ is Euler's gamma function [2], Vol. 1 and $\gamma
> \hbox{max}(0, -{\rm Re} z)$. Reciprocally, we have the direct
Mellin transform of the modified Bessel function, namely
$$\Gamma(s+z)\Gamma(s)= 2 \int_0^\infty K_z(2\sqrt x) x^{s+z/2-1}
dx.\eqno(1.3)$$
The left-hand side of (1.2) has the following asymptotic behavior
near the origin $x \to 0+$
\begin{displaymath}
 x^{z/2}K_z(2\sqrt x)= \begin{cases}
      O(1), & {\hbox{if}\  {\rm Re} z >0}, \\
      O(x^{{\rm Re} z}), & {\hbox{if}\  {\rm Re} z < 0}, \\
      O\left(\log\left({1\over x}\right)\right), & {\hbox{if}\  z=0}\end{cases}
\end{displaymath}
and $x^{z/2}K_z(2\sqrt x)= O(e^{-2\sqrt x} x^{({\rm Re} z- 1/2)/2}),
\ x \to +\infty$.

Let us consider the following integral transformation with respect
to an index $z \in \mathbb{C}$ of the modified Bessel function
$$(Ff)(z)=  2\int_0^\infty  x^{z/2}K_z(2\sqrt x)f(x)dx.\eqno(1.4)$$
This transformation looks like the Kontorovich-Lebedev transform
[5], [8], [9]. However, it is a completely different operator and
cannot be reduced to the Kontorovich-Lebedev integral by any change
of variables and functions. As far as the author is aware, the
transform (1.4) was not studied yet, taking into account his mapping
properties and inversion formula in an appropriate class of
functions.

Our goal is to do this involving a special class of functions
related to the Mellin transform and its inversion, which was
introduced in [7]. Indeed, we have

{\bf Definition 1}. Denote by ${\cal M}^{-1}(L_c)$ the space of
functions $f(x), x \in \mathbb{R}_+$,  representable by inverse
Mellin transform of integrable functions $f^{*}(s) \in L_{1}(c)$ on
the vertical line $c =\{s \in \mathbb{C}: {\rm Re} s= c_0\}$:
$$ f(x) = {1\over 2\pi i} \int_c f^{*}(s)x^{-s}ds.\eqno(1.5)$$

The  space ${\cal M}^{-1}(L_c)$  with  the  usual operations  of
addition   and multiplication by scalar is a linear vector space. If
the norm in ${\cal M}^{-1}(L_c)$ is introduced by the formula
$$ \big\vert\big\vert f \big\vert\big\vert_{{\cal
M}^{-1}(L_c)}= {1\over 2\pi }\int^{+\infty}_{-\infty} |
f^{*}\left(c_0 +it\right)| dt,\eqno(1.6)$$
then it becomes  a Banach space.

 {\bf Definition 2 ([7], [8])}.  Let $c_1, c_2 \in \mathbb{R}$
be such that $2 \hbox{sign}\ c_1 + \hbox{sign}\  c_2 \ge 0$. By
${\cal M}_{c_1,c_2}^{-1}(L_c)$ we denote the space of functions
$f(x), x \in \mathbb{R}_+$, representable in the form (1.5), where
$s^{c_2}e^{\pi c_1|s|} f^*(s) \in L_1(c)$.

It is a Banach space with the norm
$$ \big\vert\big\vert f \big\vert\big\vert_{{\cal
M}_{c_1,c_2}^{-1}(L_c)}= {1\over 2\pi }\int_{c} e^{\pi c_1|s|}
|s^{c_2} f^{*}(s) ds|.$$
In particular, letting $c_1=c_2=0$ we get the space ${\cal
M}^{-1}(L_c)$. Moreover, it is easily seen the inclusion

$${\cal M}_{d_1,d_2}^{-1}(L_c) \subseteq {\cal
M}_{c_1,c_2}^{-1}(L_c)$$ when $2 \hbox{sign}(d_1- c_1) + \hbox{sign}
(d_2-c_2) \ge 0$.

\section{Mapping properties and an inversion formula}

We begin with the following result.

{\bf Theorem 1}. {\it Let $f \in {\cal M}^{-1}(L_c)$ and $c_0 < 1$.
Then transformation $(1.4)$ is well-defined and $(Ff)(z)$ is
analytic in the half-plane ${\rm Re} z > c_0 -1$. Further,
$$(Ff)(z)={1\over 2\pi i}
\int_{c_0-i\infty}^{c_0 +i\infty} \Gamma(1-s+z)\Gamma(1-s)f^*(s)
ds,\eqno(2.1)$$
and the operator  $F: {\cal M}^{-1}(L_c) \to L_1({\rm Re} z-i\infty,
{\rm Re} z +i\infty),\ {\rm Re z}
> c_0 - 1$ is bounded with the norm satisfying the estimate}
$$||F|| \le  \Gamma(1-c_0) \int_{-\infty}^\infty |\Gamma(1-c_0+ {\rm Re z}+ i\tau)|d\tau.$$

\begin{proof} In fact, substituting (1.5) into (1.4) and changing
the order of integration by Fubini's theorem, we call (1.3) to prove
(2.1). The inversion of the order of integration is guaranteed by
the estimate (see (1.3))
$$2\int_0^\infty  \left|x^{z/2}K_z(2\sqrt x)\right| \int_c |f^{*}(s)x^{-s}ds|
dx$$$$ \le 2\int_0^\infty x^{({\rm Re}z- 2c_0)/2}K_{{\rm
Re}z}(2\sqrt x) dx \int_c |f^{*}(s)ds|$$
$$= \Gamma\left(1- c_0+ {\rm Re} z\right)\Gamma(1-c_0) \int_c |f^{*}(s)ds|< +\infty,\ {\rm Re} z >
c_0- 1, \ c_0 < 1$$
and the asymptotic behavior of the modified Bessel function at
infinity and near the origin (see above). Furthermore, integral
(1.4) converges absolutely in the half-plane ${\rm Re} z > c_0 - 1$
and uniformly in ${\rm Re} z \ge a_0 > c_0 1$. Since for each $x >0$
the function $  x^{z/2}K_z(2\sqrt x)$ is analytic by $z$, we have
that $F(z)$ is well-defined and represents an analytic function in
the half-plane ${\rm Re} z > c_0 - 1$. Finally, the straightforward
estimate takes place
$$||Ff||_1= \int_{-\infty}^{\infty} |(Ff)({\rm Re} z + i\tau)|d\tau $$
$$\le {1\over 2\pi} \int_{-\infty}^{\infty}
\int_{-\infty}^\infty |\Gamma(1- c_0+ {\rm Re
z}+i(\tau-t))\Gamma(1-c_0-it)f^*(c_0+it)| dtd\tau$$
$$\le \Gamma(1-c_0)\ ||f||_{{\cal M}^{-1}(L_c)}
\int_{-\infty}^\infty |\Gamma(1- c_0+ {\rm Re z}+ i\tau)|d\tau,$$
which completes the proof of the theorem.
\end{proof}

For the subspace $ {\cal M}_{0,n}^{-1}(L_c) \subseteq {\cal
M}^{-1}(L_c), \ n \in \mathbb{N}_0$ we have

{\bf Theorem 2}. {\it Let $n \in \mathbb{N}_0, \ f \in {\cal
M}_{0,n}^{-1}(L_c)$ and $c_0 < 1-n$. Then $f(x), \ x \in
\mathbb{R}_+$  is $n$ times continuously differentiable,
$\left(Ff^{(n)}\right)(z)$ is analytic in the half-plane ${\rm Re} z
> c_0 +n -1$ and $\left(Ff^{(n)}\right)(z)= (Ff)(z-n)$. Finally, for
any arbitrary $y \in \mathbb{R}_+$ the following representation
holds
$$(Ff)_y(z)= 2\int_y^\infty  x^{z/2}K_z(2\sqrt x)f(x)dx$$
$$= 2\sum_{m=0}^{n-1} (-1)^{m}  y^{(z+m+1)/2}K_{z+m+1}(2\sqrt y)
f^{(m)}(y) + (-1)^n \left(F f^{(n)}\right)_y(z+n), \ y
>0,\eqno(2.2)$$

where the empty sum ($n=0$) is equal to zero.}

\begin{proof} Clearly, from representation (1.1) after
differentiation and integration $n$ times with respect to $x$ under
the integral sign we come out, accordingly, with the identities
$$2{d^n\over dx^n}\left[x^{z/2} K_z(2\sqrt x)\right] = (-1)^n \int_0^\infty e^{-t - {x\over
t}}t^{z-n -1} dt = 2(-1)^n x^{(z-n)/2} K_{z-n}(2\sqrt
x),\eqno(2.3)$$
$${2\over (n-1)!}\int_y^\infty (x-y)^{n-1} x^{z/2}K_z(2\sqrt x) dx =  \int_0^\infty e^{-t - {y\over
t}}t^{z+n -1} dt$$
$$= 2y^{(z+n)/2}K_{z+n}(2\sqrt y) , \ y > 0.\eqno(2.4)$$
Further, from Definition 2 it follows that $f$ is $n$ times
continuously differentiable and via (1.5) it has
$$f^{(n)}(x)= {(-1)^n \over 2\pi }\int_{c}  (s)_n f^{*}(s) x^{-s-n}
ds,\eqno(2.5)$$
where $(a)_n$ is Pochhammer's symbol. Hence, considering $(F
f^{n})(z)$, we integrate by parts in the corresponding integral
(1.4), taking into account that the integrated terms are vanished
owing to the asymptotic behavior of the modified Bessel function,
the estimate $f^{(n)}= O(x^{-c_0-n}), \ x >0$ (see (2.5)) and limit
relations
$$\lim_{x \to 0+}x^{1-c_0-j ({\rm Re} z-  i)/2} K_{{\rm Re} z-i}(2\sqrt
x)=0, \ i, j \in \mathbb{N}_0, i+j= n,$$
which take place by virtue of the conditions $c_0 < 1-n, \ {\rm Re}
z > c_0 +n -1$. Thus calling (2.3) we prove the equality
$\left(Ff^{(n)}\right)(z)= (Ff)(z-n)$ and similar to the proof of
Theorem 1 we easily justify the  analyticity of $G(z)=
\left(Ff^{(n)}\right)(z)$ in the half-plane ${\rm Re} z
> c_0 +n -1$. Finally, the proof of (2.2) follows immediately, appealing  to (2.4) and
integrating $n$ times by parts in its left-hand side.
\end{proof}

In order to establish an inversion formula for the transformation
(1.4) we employ an operational technique, which was used formally by
Sneddon [5], Ch. 6 to deduce the inversion formula for the
Kontorovich-Lebedev transform. We start multiplying both sides of
the equality (2.1) by $x^z,\ x>0$ and integrating with respect to
$z$ over the line $(\gamma-i\infty, \gamma +i\infty), \ \gamma >
c_0-1$. Changing the order of integration in the right-hand side of
the obtained equality, which is possible via Theorem 1 and
calculating the corresponding inverse Mellin transform of the
gamma-function, we derive
$$\int_{\gamma-i\infty}^{\gamma +i\infty} (Ff)(z) x^z dz ={1\over 2\pi i}
\int_{c_0-i\infty}^{c_0 +i\infty} \int_{\gamma-i\infty}^{\gamma
+i\infty} \Gamma(1-s+z)\Gamma(1-s)f^*(s) x^z dz ds$$
$$= e^{-1/x} \int_{c_0-i\infty}^{c_0 +i\infty} \Gamma(1-s)f^*(s) x^{s-1}ds.$$
Hence, taking into account that $f \in  {\cal M}^{-1}(L_c)$, we
apply the Mellin -Parseval identity to the right-hand side of the
latter equality. Thus
$${1\over 2\pi i}\int_{\gamma-i\infty}^{\gamma +i\infty} (Ff)(z)  e^{1/x}  x^z dz
= (Lf)(x)=  \int_0^\infty e^{-xt}f(t) dt,\ x >0\eqno(2.6)$$
and the right-hand side of the latter equality represents the
Laplace transform denoted by $(Lf)(x)$. In the meantime, relation
(2.15.5.4) in [4], Vol. 2 gives the key integral involving the
modified Bessel function of the third kind $I_\nu(w)$ [2], Vol. II
$$e^{1/x} x^z= \int_0^\infty e^{-xt} I_{-(1+z)}\left(2\sqrt t\right)\  t^{-(1+z)/2}dt,\  x >0, \ {\rm Re} z < 0.$$
Substituting this integral into the left-hand side of (2.6) and
assuming an additional condition
$$(Ff)(\gamma +i\tau) \in L_1\left(|\tau|> 1; |\tau|^{\gamma+1/2}\
e^{\pi|\tau|/2}d\tau\right),\ \gamma \in (c_0-1, 0),\eqno(2.7)$$
we change of integration by Fubini's theorem and arrive at the
equality
$$\int_0^\infty e^{-xt} {1\over 2\pi i} \int_{\gamma-i\infty}^{\gamma +i\infty} I_{-(1+z)}\left(2\sqrt t\right)\  t^{-(1+z)/2}\ (Ff)(z) dz dt =
\int_0^\infty e^{-xt}f(t) dt,\ x > 0.\eqno(2.8)$$
Indeed, the motivation of the inversion of the order of integration
in (2.8)  is given due to the representation of the modified Bessel
function $I_{-(1+z)}\left(2\sqrt t\right)$ in terms of the series
$$I_{-(1+z)}\left(2\sqrt t\right)= \sum_{n=0}^\infty {t^{n-
(1+z)/2}\over n!\ \Gamma(n-z) }\eqno(2.9)$$
and an absolute integrability by $\tau \in \mathbb{R}$ of the
product $(Ff)(\gamma+i\tau)I_{-(1+\gamma+i\tau)}\left(2\sqrt
t\right)$ under condition (2.7), since $\Gamma(n-\gamma-i\tau)=
O(|\tau|^{n-\gamma-1/2} e^{-\pi|\tau|/2}), \ |\tau| \to \infty,\ n
\in \mathbb{N}_0$ via Stirling's formula [2], Vol. I.  Finally, we
observe that equality (2.8) is true for all $x >0$, where functions
under the convergent Laplace integrals in its both sides are
continuous on $\mathbb{R}_+$ owing to  condition $f \in {\cal
M}^{-1}(L_c)$ and assumption (2.7). Therefore one can cancel the
Laplace transform in (2.8) by virtue of the uniqueness theorem (see
in [3]) to get the inversion formula for the Kontorovich-Lebedev
transformation (1.4). Thus we have proved

{\bf Theorem 3.}  {\it Let $f(t) \in {\cal M}^{-1}(L_c), \ c_0 < 1$
and condition $(2.7)$ holds. Then for all $t >0$ the following
inversion formula for the transformation $(1.4)$ takes place
$$f(t)= {1\over 2\pi i} \int_{\gamma-i\infty}^{\gamma +i\infty} I_{-(1+z)}\left(2\sqrt t\right)  t^{-(1+z)/2}\ (Ff)(z)
dz,\ \gamma \in (c_0-1, 0),\eqno(2.10)$$
where the integral is absolutely convergent.}

\section{Expansion of an arbitrary function in terms of the Kontorovich-Lebedev-like integral}

In this section we will prove that any function from the  space
${\cal M}_{0, (|\varepsilon|+\varepsilon)/2}^{-1}(L_c), \ c_0<1, \
2c_0-1 < \varepsilon < c_0$ can be expanded in terms of the
following integral
$$f(x)= {1\over \pi i} {d\over dx}\int_{\gamma-i\infty}^{\gamma +i\infty} I_{-z}\left(2\sqrt x\right)  x^{-z/2}\
\int_0^\infty  t^{z/2}K_z(2\sqrt t)f(t)\ dt dz,\ x >0,\eqno(3.1)$$
where $\gamma$ is taken from the interval $(c_0-1, (\varepsilon
-1)/2)$.

Precisely, we have

{\bf Theorem 4}. {\it Let $c_0 <1, \ 2c_0-1 < \varepsilon < c_0$ and
$f \in {\cal M}_{0,(|\varepsilon|+\varepsilon)/2}^{-1}(L_c)$. Then
for any $x >0$ formula $(3.1)$ is true, where the interior integral
with respect to $t$ converges absolutely and the exterior integral
by $z$ is understood in the improper sense of Riemann}.

\begin{proof} In fact, since ${\cal M}_{0,(|\varepsilon|+\varepsilon)/2}^{-1}(L_c) \subseteq {\cal
M}^{-1}(L_c)$, the absolute convergence of the interior integral in
(3.1)  follows from Theorem 1. Moreover, equality (2.1) holds. Hence
writing the modified Bessel function $I_{-z}\left(2\sqrt x\right)$
similar to (2.9) and substituting the right-hand side of (2.1) into
(3.1), we come out with the following iterated integral
$$I(x)= - {1\over 4\pi^2}\int_{\gamma-i\infty}^{\gamma +i\infty} \sum_{n=0}^\infty {x^{n-
z}\over n!\ \Gamma(1+n-z) }\ \int_{c_0-i\infty}^{c_0 +i\infty}
\Gamma(1-s+z)\Gamma(1-s)f^*(s) ds dz.\eqno(3.2)$$
Meanwhile, appealing to the Stirling formula for gamma-functions
[2], Vol. I, we find for any $n \in \mathbb{N}$
$$\left|\frac{\Gamma(1-s+z)}{\Gamma(1+n-z) }\right|=
\left|B(1-s+z, s-\varepsilon) B(n, 1-z)\frac{\Gamma(1-\varepsilon
+z)}{\Gamma(s-\varepsilon)\Gamma(1-z)(n-1)! }\right|$$
$$\le {B(1-c_0+\gamma, c_0-\varepsilon)\Gamma(1-\gamma)\over \Gamma(1+n-\gamma)}\left|\frac{\Gamma(1-\varepsilon
+z)}{\Gamma(s-\varepsilon)\Gamma(1-z)}\right|= O\left({|z|^{2\gamma-
\varepsilon}\over |\Gamma(s-\varepsilon)|}\right),\ \ |{\rm Im} z|
\to \infty,$$
where $B(a,b)$ is Euler's beta-function,\ $c_0-1 < \gamma <
(\varepsilon-1)/2$. Hence from (3.2) for each fixed $x >0$ we obtain
the estimate
$$\int_{\gamma-i\infty}^{\gamma +i\infty} \sum_{n=0}^\infty \left|{x^{n-
z}\over n!\ \Gamma(1+n-z) }\right|\ \int_{c_0-i\infty}^{c_0
+i\infty} \left|\Gamma(1-s+z)\Gamma(1-s)f^*(s) ds dz\right|$$
$$\le  \ B(1-c_0+\gamma, c_0-\varepsilon)\Gamma(1-\gamma) x^{-\gamma/2} I_{-\gamma}(2\sqrt
x)$$
$$\times  \int_{\gamma-i\infty}^{\gamma +i\infty} \left|\frac{\Gamma(1-\varepsilon
+z)}{\Gamma(1-z)}\right|\ \int_{c_0-i\infty}^{c_0 +i\infty}
\left|{\Gamma(1-s)\over \Gamma(s-\varepsilon)} f^*(s) ds dz\right|$$
$$= O\left(\int_{\gamma-i\infty}^{\gamma +i\infty} |z|^{2\gamma-
\varepsilon}|dz| \int_{c_0-i\infty}^{c_0 +i\infty} |s|^\varepsilon |
f^*(s) ds| \right)< +\infty.$$
Consequently, the change of the order of integration and summation
is possible in (3.2). After calculation of the integral with respect
to $z$ using relation (8.4.19.1) in [4], Vol. 3 it becomes
$$I(x)= {1\over 2\pi i} \int_{c_0-i\infty}^{c_0 +i\infty} \sum_{n=0}^\infty {x^{n}\over n!}J_{n+1-s}(2\sqrt x)
\ \Gamma(1-s)f^*(s)x^{(1-s)/2} ds,\eqno(3.3)$$
where $J_\mu(w)$ is the Bessel function of the first kind [2], Vol.
II. But the series inside (3.3) is calculated in [4], Vol. 2,
relation (5.7.6.7), namely
$$\sum_{n=0}^\infty {x^{n}\over n!}J_{n+1-s}(2\sqrt x)=
{x^{(1-s)/2}\over \Gamma(2-s)}.$$
Thus substituting this value into (3.3) and applying the reduction
formula for gamma-function, we arrive at the equality
$$I(x)= {1\over 2\pi i} \int_{c_0-i\infty}^{c_0 +i\infty} f^*(s){x^{1-s}\over 1-s} ds.\eqno(3.4)$$
Hence the differentiation with respect to $x >0$ under integral sign
in (3.4) is permitted via the absolute and uniform convergence since
$f^*(s) \in L_1(c)$ (see Definition 2). Thus we establish equality
(3.1) and complete the proof.
\end{proof}

As we see, expansion (3.1) generates the following reciprocal
inversion formula of the index transform (1.4)
$$f(x)= {1\over 2\pi i} {d\over dx}\int_{\gamma-i\infty}^{\gamma +i\infty}
I_{-z}\left(2\sqrt x\right) x^{-z/2} (Ff)(z) dz, \ x >0.\eqno(3.5)$$

{\bf Corollary 1.} {\it Let, in addition, condition $(2.7)$ hold.
Then formula $(3.5)$ can be written in the form $(2.10)$.}

\begin{proof} Indeed, in this case the differentiation under
integral sign in (3.5) is allowed via the absolute and uniform
convergence. Hence using the identity for derivatives of Bessel
functions [2], Vol. II
$${d\over dx} \left[I_{-z}\left(2\sqrt x\right)  x^{-z/2}\right]=
I_{-(z+1)}\left(2\sqrt x\right)  x^{-(z+1)/2},$$
we arrive at the result.
\end{proof}

{\bf Corollary 2}. {\it Let  $c_0 <1, \ 2c_0-1 < \varepsilon < c_0$
and $f \in {\cal M}_{0,(|\varepsilon|+\varepsilon)/2}^{-1}(L_c)$.
Then the homogeneous integral equation
$$\int_0^\infty  t^{z/2}K_z(2\sqrt t)f(t) dt=0$$
has only the trivial solution.}

Expansion (3.1) gives a new source of index integrals involving the
modified Bessel function $I_\nu(w)$. It can be obtained employing
the corresponding integrals (1.4) for concrete functions $f$ from
[4], Vol. 2. In fact,  making a simple substitution in (1.4) and
then using relation (2.16.6.4) in [4], Vol. 2 we calculate the value
of the index integral
$${1\over 2\pi i} \int_{\nu-i\infty}^{\nu +i\infty}
I_{-z}\left(2\sqrt x\right){\Gamma(z)\over 2z+1} x^{-z/2}\ dz=
e^{-2\sqrt x},\ x >0; \ \nu < 1/2.$$
Meanwhile, relation (2.16.33.2) in [4], Vol. 2 leads us to the value
of the reciprocal index integral
$${1\over 4\pi i} \int_{\nu-i\infty}^{\nu +i\infty}
I_{-z}\left(2\sqrt x\right){\Gamma\left(z+ {\mu\over
2}\right)\Gamma\left(z - {\mu\over 2}\right)\over \Gamma(z+1)}
x^{-z/2}\ dz$$ $$=  K_\mu(2\sqrt x)\left[\Gamma\left(1+ {\mu\over
2}\right)\Gamma\left(1 - {\mu\over 2}\right)\right]^{-1}, \ x >0;\
|{\rm Re} \mu| /2 < \nu < 1/2.$$
More curious example can be calculated, for instance,  via relation
(2.16.8.4) in [4], Vol. 2.  Indeed, we have
$${1\over 2\pi i} \int_{\nu-i\infty}^{\nu +i\infty}
I_{-z}\left(2\sqrt x\right)W_{-z/2, \ (z-1)/2}\left({1\over
4p}\right)\Gamma\left(z\right) (4px)^{-z/2}\ dz$$
$$=  e^{-4px- {1\over 8p}}, \ x, p >0;\  1/2 < \nu < \varepsilon +
1/2, \ \varepsilon \in (0, 1/4).$$
where $W_{\mu,\nu}(w)$ is the Whittaker function [2], Vol. I.

\section{A convolution operator and integral equations of the convolution type}

In this section we will construct a convolution operator, which is
related to the transformation (1.4) and the Mellin transform [6]
$$\left({\cal M} f\right)(z)=\int_0^\infty f(x)x^{z-1}dx.\eqno(4.1)$$
Our construction will be based on the convolution properties of the
Mellin transform in $L_1$ (see [5], [6], Th. 44)  and representation
(1.1). Indeed, considering (1.1) of the same parameter $z$ and
different positive arguments $x$ and $y$, we deduce the following
representation of the product of these integrals, namely
$$4(xy)^{z/2}K_z(2\sqrt x)K_z(2\sqrt y)= \int_0^\infty e^{-t - {x\over
t}}t^{z-1} dt\int_0^\infty e^{-u - {y\over u}}u^{z-1} du$$
$$=\int_0^\infty v^{z-1} \left(\int_0^\infty e^{-\ {t( y+v)\over v} -\  {x+v\over
t}}{dt\over t}\right) dv= 2\int_0^\infty
K_0\left(2\sqrt{{(x+v)(y+v)\over v}}\right) v^{z-1} dv,$$
 where the change of the order of integration is allowed by the
 Fubini theorem via the absolute convergence. So we find the product
 integral formula for the kernel of transformation (1.4)
$$2(xy)^{z/2}K_z(2\sqrt x)K_z(2\sqrt y)= \int_0^\infty
K_0\left(2\sqrt{{(x+v)(y+v)\over v}}\right) v^{z-1} dv,\ (x,y) \in
\mathbb{R}_+^2, \ z \in \mathbb{C}.\eqno(4.2)$$

{\bf Definition 3}. We will call the following bilinear form
$(f*g)(x), \ x \in \mathbb{R}_+$
$$(f*g)(x)= 2 \int_{\mathbb{R}_+^2}
K_0\left(2\sqrt{{(x+u)(x+v)\over x}}\right)f(u)g(v)dudv\eqno(4.3)$$
a convolution operator for the transformation (1.4) whenever it
exists.

Let us consider the weighted $L_1$-space  $L_1(\mathbb{R}_+; 2
x^{\alpha/2}K_\alpha(2\sqrt x) dx), \ \alpha \in \mathbb{R}$ with
the norm
$$||f||_{L_1(\mathbb{R}_+; 2
x^{\alpha/2}K_\alpha(2\sqrt x) dx)}= 2 \int_0^\infty |f(x)|
x^{\alpha/2}K_\alpha(2\sqrt x)dx.$$

Similar to (2.6), we prove first the composition representation of
the transformation (1.4) in terms of the Mellin and Laplace
integrals.

{\bf Theorem 5}. {\it Let $f \in L_1(\mathbb{R}_+;
x^{(\alpha-|\alpha|)/2} dx), \ \alpha \neq 0$. Then $(Ff)(z)$ is
analytic in the right half-plane
\begin{displaymath}
{\rm Re}\ z \ge \begin{cases}
      0, & {\hbox{if}\  \alpha > 0}, \\
      \alpha, & {\hbox{if}\  \alpha < 0} \end{cases}
\end{displaymath}
and  can be represented there by the composition of the Mellin and
Laplace transforms as follows}
$$(Ff)(z)= {\cal M}\circ \left(e^{-t} (Lf)(1/t)\right) (z).\eqno(4.4)$$

\begin{proof} The proof is straightforward by Fubini's theorem with
the use of integral representation (1.1), asymptotic behavior of the
modified Bessel function  and the estimates
$$|x^{z/2}K_z(2\sqrt x)| \le x^{{\rm Re} z/2}K_{{\rm Re} z}(2\sqrt
x)\le C x^{\beta/2}K_{\beta}(2\sqrt x),\ x >0,\eqno(4.5)$$
where $C >0$ is an absolute constant when
\begin{displaymath}
{\rm Re}\ z \ge \begin{cases}
      0, & {\hbox{if}\  \beta \ge 0}, \\
      \beta, & {\hbox{if}\  \beta < 0}, \end{cases}
\end{displaymath}
$$2x^{\alpha/2}K_{\alpha}(2\sqrt x)\le x^{(\alpha- |\alpha|)/2} \Gamma(|\alpha|), \
\alpha \neq 0,\eqno(4.6)$$
$$\int_0^\infty  |x^{z/2}K_z(2\sqrt x)|f(x)|dx \le C \int_0^\infty  |f(x)| \int_0^\infty e^{-t - {x\over t}}t^{\alpha -1} dtdx$$
$$\le C \ \Gamma(|\alpha|) \int_0^\infty x^{(\alpha-|\alpha|)/2} |f(x)|dx < \infty.$$
\end{proof}

{\bf Theorem 6}. {\it Let $f,g \in L_1(\mathbb{R}_+; 2
x^{\alpha/2}K_\alpha(2\sqrt x) dx), \ \alpha \in \mathbb{R}$. Then
convolution $(4.3)$ exists and belongs to the space
$L_1(\mathbb{R}_+; x^{\alpha-1}dx)$, satisfying the Young type
inequality
$$||f*g||_{L_1(\mathbb{R}_+; x^{\alpha-1} dx)} \le ||f||_{L_1(\mathbb{R}_+; 2
x^{\alpha/2}K_\alpha(2\sqrt x) dx)}||g||_{L_1(\mathbb{R}_+; 2
x^{\alpha/2}K_\alpha(2\sqrt x) dx)}.\eqno(4.7)$$
Moreover, this form is commutative and the following factorization
equality holds in terms of transformations $(1.4)$, $(4.1)$

$$\left({\cal M}(f*g)\right)(z)= (Ff)(z)(Fg)(z),\eqno(4.8)$$
where $z$ belongs to the half-plane}
\begin{displaymath}
{\rm Re}\ z \ge \begin{cases}
      0, & {\hbox{if}\  \alpha \ge 0}, \\
      \alpha, & {\hbox{if}\  \alpha < 0}. \end{cases}
\end{displaymath}

\begin{proof}  In fact, the existence of the convolution (4.3) for almost all $x >0$ follows
from Fubini's theorem and the estimate
$$\int_0^\infty |(f*g)(x)| x^{\alpha- 1}dx \le 2 \int_0^\infty x^{\alpha-1}
\int_{0}^\infty \int_0^\infty  K_0\left(2\sqrt{{(x+u)(x+v)\over
x}}\right)|f(u)g(v)|dudv dx $$
$$=4\int_0^\infty  u^{\alpha/2}K_\alpha(2\sqrt u)|f(u)|du \int_0^\infty
v^{\alpha/2}K_\alpha(2\sqrt v)|g(v)|dv. $$ This also drives us to
the Young type inequality (4.6). Hence the factorization equality
(4.7) is an immediate consequence of (4.2), (4.5) with
$\beta=\alpha$  and straightforward calculations.
\end{proof}

Letting $\alpha=1$ and using inequality (4.6)  we obtain as a
corollary the $L_1$-property of the convolution (4.3).

{\bf Corollary 3}. {\it Let $f,g \in L_1(\mathbb{R}_+; dx)$. Then
convolution $(4.3)$ exists and belongs to $L_1(\mathbb{R}_+; dx)$,
yielding the corresponding Young inequality
$$||f*g||_{L_1} \le ||f||_{L_1}||g||_{L_1}.\eqno(4.9)$$
Moreover, the convolution  is commutative and associative,
satisfying the factorization equality $(4.8)$ in the half-plane
${\rm Re} z \ge 0$.}

Further, appealing to Corollary 2 we prove  an analog of
Titchmarsh's theorem about the absence of divisors of zero for
convolution (4.3).

{\bf Theorem 7.} {\it Let $f,g \in L_1(\mathbb{R}_+; dx)$.  Then the
equality $(f*g)(x)=0$ yields that at least one of the functions
$f(x)$ and $g(x)$ is equal to zero for all $x >0$.}

\begin{proof} In fact, both functions $(Ff)(z), (Fg)(z)$ are
analytic in the half plane ${\rm Re} z >0$ and via equality (4.8) at
least on of them is identically equal to zero. Then the result
follows from Theorem 5 due to the uniqueness theorems in $L_1$ for
the Mellin and Laplace transforms.
\end{proof}

The Parseval type equality for convolution (4.3) is an immediate
consequence of the Plancherel $L_2$-theory of the Mellin transform
[6].  We have

{\bf Theorem 8}. {\it Let $f,g \in L_1(\mathbb{R}_+; dx)$. Then
$(f*g)(x) \in L_2(\mathbb{R}_+; x^{2\alpha -1}dx), \ \alpha >0$ and
the Parseval type equality holds}
$$\int_0^\infty |(f*g)(x)|^2 x^{2\alpha-1} dx= {1\over
2\pi}\int_{-\infty}^\infty \left|(Ff)(\alpha+it)(Fg)(\alpha
+it)\right|^2 dt.\eqno(4.10)$$

\begin{proof} Indeed by virtue of the generalized Minkowskii
inequality and relation (8.4.23.27) in [4], Vol. 3 we derive
$$||f*g||_{L_2(\mathbb{R}_+; x^{2\alpha -1}dx)}= 2\left(\int_0^\infty
\left|\int_{\mathbb{R}_+^2} K_0\left(2\sqrt{{(x+u)(x+v)\over
x}}\right)f(u)g(v)dudv\right|^2 x^{2\alpha-1} dx\right)^{1/2}$$
$$\le 2\int_{\mathbb{R}_+^2} |f(u)g(v)|
\left(\int_0^\infty K_0^2\left(2\sqrt{{(x+u)(x+v)\over x}}\right)x^{2\alpha-1}dx\right)^{1/2}dudv$$
$$\le \left(\int_0^\infty
K_0^2\left(2\sqrt{x}\right)x^{2\alpha-1}dx\right)^{1/2}
||f||_{L_1}||g||_{L_1}=
2^{-2\alpha-1/2}\pi^{1/4}{\Gamma^{3/2}(2\alpha)\over
\Gamma^{1/2}(2\alpha+1/2)}||f||_{L_1}||g||_{L_1} < \infty.$$
Hence factorization equality (4.8) and Theorem 71 in [6] give the
result.
\end{proof}

Finally, let us consider a class of convolution integral equations
of the first kind generated by (4.3)
$$  \int_{0}^\infty k_h(x,y)f(y)dy = g(x), \ x> 0,\eqno(4.11)$$
where
$$k_h(x,y)= 2\int_0^\infty  K_0\left(2\sqrt{{(x+y)(x+u)\over
x}}\right)h(u)du,\eqno(4.12)$$
 $h, g$ are given functions and $f$ is to be determined.

{\bf Theorem 9.}  {\it Let $f \in {\cal M}^{-1}(L_c), \ c_0 < 1$, $h
\in L_1(\mathbb{R}_+; 2 x^{\alpha/2}K_\alpha(2\sqrt x) dx), \ 0
> \alpha > c_0-1$ and $g \in L_1(\mathbb{R}_+; x^{\alpha-1}dx)$. Let also  transformation
$(1.4)$ of $h \ (Fh)(z)$ has no zeros in the strip ${\rm Re} z \in
(\alpha, 0)$ and the quotient $\left({\cal M} g\right)(z)/(Fh)(z)$,
where $\left({\cal M} g\right)(z)$ is the Mellin transform $(4.1)$
of $g$,   satisfies condition $(2.7)$ in this strip.  Then a
solution of integral equation $(4.12)$ has the form}
$$f(x)= {1\over 2\pi i} \int_{\gamma-i\infty}^{\gamma +i\infty} I_{-(1+z)}\left(2\sqrt x\right)  x^{-(1+z)/2}\
 {\left({\cal M} g\right)(z)\over (Fh)(z)}dz,\ x >0, \gamma \in (\alpha, 0).\eqno(4.13)$$
\begin{proof} Clearly, by straightforward estimate of the norm we
verify that if $f \in {\cal M}^{-1}(L_c), \ c_0 < 1$ and $\alpha \in
(c_0-1,\ 0)$, then $f \in L_1(\mathbb{R}_+; 2
x^{\alpha/2}K_\alpha(2\sqrt x) dx)$. Therefore Theorem 6 and formula
(4.8) are valid for convolution $(f*h)(x)$. Hence since $(Fh)(z)\neq
0$ in the strip ${\rm Re} z \in (\alpha, 0)$ it has the equality
$$(Ff)(z)= {\left({\cal M}g\right)(z)\over (Fh)(z)}.$$
Consequently, appealing to Theorem 3 and formula (2.10), we complete
the proof of the theorem.
\end{proof}

An interesting example of the equation (4.12) and its solution can
be found, taking, for instance,  $h(x)= x^{-1/2}$. In this case one
can calculate the kernel (4.12) via relation (2.16.3.10) in [4],
Vol. 2 and we obtain

$$k_h(x,y)= {\pi \sqrt x \over \sqrt{x+y}} e^{-2\sqrt{x+y}}.$$
Moreover, it has $(Fh)(z)=\sqrt \pi \Gamma(z+1/2)$ by formula (1.3).
Hence Theorem 9 says that a solution of the integral equation
$$\pi \sqrt x  \int_{0}^\infty {e^{-2\sqrt{x+y}}\over \sqrt{x+y}}f(y)dy = g(x), \ x> 0,$$
is given by the integral
$$f(x)= {1\over 2\pi\sqrt \pi i} \int_{\gamma-i\infty}^{\gamma +i\infty} {I_{-(1+z)}\left(2\sqrt x\right)\over \Gamma(z+1/2)}  x^{-(1+z)/2}\
 \left({\cal M} g\right)(z)dz,\ x >0,$$
where $\gamma \in (\alpha, 0)$, $\alpha$  is chosen  from the
interval   $\alpha \in (\hbox{max}(c_0-1, -1/2),\  0)$ and

$$\left({\cal M} g\right)(\gamma +i\tau) \in L_1\left(|\tau|> 1;\  |\tau|^{1/2}\
e^{\pi|\tau|}d\tau\right).$$

This example can be generalized, considering $h(x)= x^{\beta-1}, \
\beta > 0$. Hence using relation (2.16.3.8) in [4], Vol. 2, we find
$$k_h(x,y)= 2\Gamma(\beta) \left({x\over \sqrt{x+y}}\right)^\beta
K_\beta\left(2\sqrt{x+y}\right).$$ Moreover,
 $(Fh)(z)=\Gamma(\beta)\Gamma(\beta+z)$ and a solution of the
 equation

$$2\Gamma(\beta) \int_{0}^\infty \left({x\over \sqrt{x+y}}\right)^\beta
K_\beta\left(2\sqrt{x+y}\right)f(y)dy = g(x), \ x> 0,$$
is
$$f(x)= {1\over 2\pi\Gamma(\beta) i} \int_{\gamma-i\infty}^{\gamma +i\infty} {I_{-(1+z)}\left(2\sqrt x\right)\over \Gamma(z+\beta)}  x^{-(1+z)/2}\
 \left({\cal M} g\right)(z)dz,\ x >0,$$
where $\gamma \in (\alpha, 0)$, $\alpha$  is chosen  from the
interval   $\alpha \in (\hbox{max}(c_0-1, -\beta),\  0)$ and

$$\left({\cal M} g\right)(\gamma +i\tau) \in L_1\left(|\tau|> 1;\  |\tau|^{1-\beta}\
e^{\pi|\tau|}d\tau\right).$$

Finally we write this solution in terms of the Neumann type series.
In fact, substituting the value of the modified Bessel function
$I_{-(1+z)}\left(2\sqrt x\right)$ by series (2.9),  we change the
order of summation and integration via the absolute convergence to
obtain

$$f(x)= {1\over \Gamma(\beta)}\sum_{n=0}^\infty {x^{n- 1}\over n!\ \Gamma(n+\beta)} \left\{x^\beta (1+x)^{-\beta-n}\right\}^{-1}g,$$
where by the symbol

$$\left\{x^\beta (1+x)^{-\beta-n}\right\}^{-1}g =
{1\over 2\pi  i} \int_{\gamma-i\infty}^{\gamma +i\infty}
{\left({\cal M} g\right)(z)\over \Gamma(n-z) \Gamma(z+\beta)}
x^{-z}dz$$ the generalized inverse Stieltjes transform is denoted
(see details in [1]).

\bigskip
\centerline{{\bf Acknowledgments}}
\bigskip
The present investigation was supported, in part,  by the "Centro de
Matem{\'a}tica" of the University of Porto.

\bigskip
\centerline{{\bf References}}
\bigskip
\baselineskip=12pt
\medskip
\begin{enumerate}

\item[{\bf 1.}\ ] Yu.A. Brychkov, Kh.-Yu. Glaeske and  O.I. Marichev, Factorization of integral
transformations of convolution type,  {\it Mathematical analysis},
{\bf 21} (1983), 3-41 (Russian).

\item[{\bf 2.}\ ]
 A. Erd\'elyi, W. Magnus, F. Oberhettinger and F.G. Tricomi,
{\it Higher Transcendental Functions}, Vols. I and  II, McGraw-Hill,
New York, London and Toronto (1953).

\item[{\bf 3.}\ ] V.A. Ditkin and A.P. Prudnikov, {\it Operational
Calculus}. Nauka, Moscow, 1975 (in Russian).

\item[{\bf 4.}\ ]  A.P. Prudnikov, Yu. A. Brychkov and O. I.
Marichev, {\it Integrals and Series: Vol. 2: Special Functions},
Gordon and Breach, New York  (1986); {\it Vol. 3: More Special
Functions}, Gordon and Breach, New York  (1990).

\item[{\bf 5.}\ ]I.N. Sneddon, {\it The Use of Integral
Transforms}, McGray Hill, New York (1972).

\item[{\bf 6.}\ ]  E.C. Titchmarsh, {\it  An Introduction to the
Theory of Fourier Integrals}, Clarendon Press, Oxford ( 1937).

\item[{\bf 7.}\ ] Vu Kim Tuan, O.I. Marichev and S.B. Yakubovich, Composition structure of integral
transformations,  {\it J. Soviet Math.}, {\bf 33} (1986), 166-169.

\item[{\bf 8.}\ ] S. B. Yakubovich and  Yu. F. Luchko, {\it The
Hypergeometric Approach to Integral Transforms and Convolutions}.
Mathematics and its Applications, 287. Kluwer Academic Publishers
Group, Dordrecht (1994).

\item[{\bf 9.}\ ] S.B. Yakubovich, {\it Index Transforms},  World
Scientific Publishing Company, Singapore, New Jersey, London and
Hong Kong (1996).

\end{enumerate}

\vspace{5mm}

\noindent S.Yakubovich\\
Department of Mathematics,\\
Faculty of Sciences,\\
University of Porto,\\
Campo Alegre st., 687\\
4169-007 Porto\\
Portugal\\
E-Mail: syakubov@fc.up.pt\\

\end{document}